\documentclass[11pt,amsfonts]{article}
\usepackage{mathrsfs}
\usepackage{graphicx}
\usepackage{latexsym}
\usepackage{amssymb}
\usepackage{amsmath}
\usepackage{layout}
\newtheorem{prop}{Proposition}

\newtheorem{definition}{Definition}

\newtheorem{remark}{Remark}

\def\real{{\mathord{{\rm J\kern-2.8pt R}}}}        
\def\inte{{\mathord{{\rm J\kern-2.8pt N}}}}

\def\sZZ{{\rm Z\kern-2.8ptem{}Z}}

\def\z{{\mathchoice
  {\sZZ}
  {\sZZ}
  {\rm Z\kern-0.30em{}Z}
  {\rm Z\kern-0.25em{}Z} }}
\def\sQQ{{\kern 0.27em \vrule height1.45ex width0.03em depth0em
          \kern-0.30em \rm Q}}
\def\qu{{\mathchoice
    {\sQQ}
    {\sQQ}
  {\kern 0.225em \vrule height1.05ex width0.025em depth0em \kern-0.25em \rm Q}
  {\kern 0.180em \vrule height0.78ex width0.020em depth0em \kern-0.20em \rm Q}
        }}
\def\sCC{{\kern 0.27em \vrule height1.45ex width0.03em depth0em
          \kern-0.30em \rm C}}
\def\complex{{\mathchoice
    {\sCC}
    {\sCC}
  {\kern 0.225em \vrule height1.05ex width0.025em depth0em \kern-0.25em \rm C}
  {\kern 0.180em \vrule height0.78ex width0.020em depth0em \kern-0.20em \rm C}
        }}

\newcommand{\ba}{\begin{array}}
\newcommand{\ea}{\end{array}}
\newcommand{\be}{\begin{equation}}
\newcommand{\ee}{\end{equation}}
\newcommand{\bea}{\begin{eqnarray}}
\newcommand{\eea}{\end{eqnarray}}
\newcommand{\beaa}{\begin{eqnarray*}}
\newcommand{\eeaa}{\end{eqnarray*}}

\def\z{\zeta}

\font\tenmath=msbm10 \font\sevenmath=msbm7 \font\fivemath=msbm5
\newfam\mathfam \textfont\mathfam=\tenmath
\scriptfont\mathfam=\sevenmath \scriptscriptfont\mathfam=\fivemath

\def \={{\buildrel {\rm (law)} \over =}}

\def\cF{{\cal F}}

\def\cH{{\cal H}}

\def\qed{ \hfill \vrule width.25cm height.25cm depth0cm\smallskip}

\newcommand{\basa}{\begin{assumption}}
\newcommand{\easa}{\end{assumption}}

\newcommand{\bas}{\begin{assum}}
\newcommand{\eas}{\end{assum}}

\newcommand{\ignore}[1]{}
\begin{document}

\renewcommand{\thefootnote}{\fnsymbol{footnote}}

\title{Wiener integrals with respect to the Hermite random field and applications to the wave equation}
\author{Jorge Clarke De la Cerda $^{1,2,3}$ \thanks{Partially supported by BECAS CHILE 2011 and the CONICYT-ECOS program C10E03}. \quad Ciprian A. Tudor
$^{3,4}$ \footnote{Supported by the CNCS grant PN-II-ID-PCCE-2011-2-0015. Associate member of the team Samm, Universit\'e de Panth\'eon-Sorbonne Paris 1.  Supported by the ANR grant "Masterie" BLAN 012103 is also acknowledged.}\vspace*{0.1in}\\ 
$^{1}$  Departamento de Matem\'atica, Facultad de Ciencias, Universidad del Bio-Bio, \\ Concepci\'on, Chile. \\
$^{2}$ Departamento de Ingenier\'ia Matem\'atica, Universidad de Concepci\'on,\\
Casilla 160-C, Concepci\'on, Chile. \\
jclarke@udec.cl \vspace*{0.1in} \\
$^{3}$  Laboratoire Paul Painlev\'e, Universit\'e de Lille 1\\
 F-59655 Villeneuve d'Ascq, France. \vspace*{0.1in} \\
 $^{4}$ 
Academy of Economical Studies, Bucharest\\
Piata Romana, nr. 6, Bucharest, Romania.\\
\quad tudor@math.univ-lille1.fr\vspace*{0.1in}}
\date{}
\maketitle

\begin{abstract}

The Hermite random field has been introduced as a limit of some weighted Hermite variations of the fractional Brownian sheet. In this work we define it as a multiple integral with respect to the standard Brownian sheet and introduce Wiener integrals with respect to it. As an application we study the wave equation driven by the Hermite sheet. We  prove the existence of the solution and we study the regularity of its sample paths, the existence of the density and of  its local times.

\end{abstract}

\vskip0.3cm

{\bf  2000 AMS Classification Numbers: }60F05, 60H05, 60G18.

 \vskip0.2cm


\section{Introduction}

\hspace*{1.1cm} The random fields or multiparameter stochastic  processes have focused a significant amount of attention among scientists due to the wide range of applications that they have. Particularly, self-similar random fields find some of their  applications in various kind of phenomena, going from hydrology and surface modeling to network traffic analysis and mathematical finance, to name a few. From other side, this type  of processes are also quite interesting when they appear as solutions to Stochastic  Partial Differential Equations (SPDE's) in several dimensions, such as the wave or heat equations.\\
\hspace*{1.1cm} A class of processes that lies in the family described above are the Hermite random fields or Hermite sheets (from now on). Inside this class we can find the well-studied fractional Brownian sheet and the Rosenblatt processes, among others.\\

The Hermite processes  of order $q\geq 1$ are  self-similar with stationary increments and live in the $q$th Wiener chaos, that is, they can be expressed as $q$ times iterated integrals with respect to the Wiener process. The class of Hermite processes includes the fractional Brownian motion (fBm) which is the only Gaussian process in this family.
Their practical aspects are striking: they provide a wide class of processes that allow to model long memory, self-similarity and H\"{o}lder-regularity,
enabling a significant deviation from the fBm and other Gaussian processes.  Since they
are non-Gaussian and self-similar with stationary increments, the Hermite processes can also be
an input in models where self-similarity is observed in empirical data which appears to be
non-Gaussian.

The Hermite sheet of order $q$ is only known in his representation as a non-central limit of a particularly normalized Hermite variations of the fractional Brownian sheet, see \cite{RST10} for the two-parameter case and \cite{Breton11} for the general $d$-parametric case. In both cases the authors also prove self-similarity, stationary increments and H\"older continuity.\\
\hspace*{1.1cm} In the present work we deal directly with the multi-parametric case building the Hermite sheet as a natural extension of the expression for the Hermite process studied as a non-central limit in \cite{DoMa} and \cite{Taqqu79}.

Fix $d\in \mathbb{N} \backslash \left\lbrace0\right\rbrace$, define $\mathbf{t}=(t_{1},t_{2}, \ldots ,t_{d}) \in \mathbb{R}^{d}$ and let $\mathbf{H}=(H_{1},H_{2}, \ldots ,H_{d}) \in (\frac{1}{2},1)^{d}$ a Hurst multi-index

\begin{eqnarray}
\label{hermite-sheet}
\nonumber
 Z^{q}_{\mathbf{H}}(\mathbf{t}) &=& c(\mathbf{H},q) \int_{\mathbb{R}^{d\cdot q}}
 \int_{0}^{t_{1}}
 \ldots \int_{0}^{t_{d}} \left( \prod _{j=1}^{q} (s_{1}-y_{1,j})_{+} ^{-\left( \frac{1}{2} + \frac{1-H_{1}}{q} \right) }
 \ldots (s_{d}-y_{d,j})_{+} ^{-\left( \frac{1}{2} + \frac{1-H_{d}}{q} \right) }   \right)  \\
\nonumber
& &  ds_{d} \ldots ds_{1} \quad
dW(y_{1,1},\ldots ,y_{d,1} )\ldots dW(y_{1,q},\ldots ,y_{d,q}) \\
&=&
c(\mathbf{H},q) \int_{\mathbb{R}^{d\cdot q}}
\int_{0}^{\mathbf{t}} \quad \prod _{j=1}^{q} (\mathbf{s}-\mathbf{y}_{j})_{+} ^{-\mathbf{\left( \frac{1}{2} + \frac{1-H}{q} \right)} }  d\mathbf{s} \quad dW(\mathbf{y}_{1})\ldots dW(\mathbf{y}_{q}),
\end{eqnarray}
the bold characters are for multidimensional quantities as indicated below in Section 2.

The above integrals are Wiener-It\^o multiple integrals of order $q$ with respect to the $d$-parametric standard Brownian sheet $(W(\mathbf{y}))_{\mathbf{y}\in \mathbb{R}^{d}}$ (see \cite{N} for the definition) and $c(\mathbf{H},q)$ is a positive
normalization constant depending only on $\mathbf{H}$ and $q$. We designate the process $\left( Z^{q}_{\mathbf{H}}(\mathbf{t}) \right)_{\mathbf{t} \in \mathbb{R}^{d} }$ as the {\em Hermite sheet} or {\em Hermite random field}.

From expression (\ref{hermite-sheet}) it is possible to note that for $d=1$ we recover the {\em Hermite process} which represents a family that has been recently studied in \cite{CTV11-Hp}, \cite{MT07} and \cite{PT10-Hp}. As a particular case ($q=1$) we recover the most known element of this family,  the {\em fractional Brownian motion}, which has been largely studied due to its various applications.  Recently, a rich theory of stochastic integration with respect to this process has been introduced  and stochastic differential equations driven by the
fractional Brownian motion have been considered for several purposes. The process obtained in (\ref{hermite-sheet-1}) for $d=1, q=2$ is known as the {\em Rosenblatt process}, it was introduced by Rosenblatt in \cite{Rose} and it has been called in this way by Taqqu in \cite{Taqqu75}. Lately, this process has been increasingly studied by his different interesting aspects like  wavelet type expansion or extremal properties, parameter estimations, discrete approximations and others potential applications (see  \cite{Abry1}, \cite{Abry2}, \cite{BT10}, \cite{CTV09},  \cite{Tudor08}). \\

As far as we know, the only well-known multiparameter process that can be obtained from (\ref{hermite-sheet}) is the {\em fractional Brownian sheet} ($d>1$ and $q=1$). This process has been recently studied as a driving noise for stochastic differential equations and  stochastic calculus with respect to it has been developed. We refer  to   \cite{ALP2002}, \cite{KP2009}, \cite{TV1} for only a few works on various aspects of the fractional Brownian sheet.

In one hand the purpose of this article is to study the basic properties of the multiparameter Hermite process  and then to introduce Wiener integrals with respect to the Hermite sheet in order to generalize and continue the line introduced in \cite{MT07} putting a new brick in the construction of stochastic calculus driven by this class of processes in several dimensions. As in \cite{Breton11} the covariance structure of the Hermite sheet is like the one of the fractional Brownian sheet, enabling the use of the same classes of deterministic integrands as in the fractional Brownian sheet profiting its  well-known properties.

Also in the aim of this work lives the idea of making an approach to the study of stochastic partial  differential equations in several dimensions driven by non-Gaussian noises, giving a specific expression for the driving noise allowing to use in a better way the properties of the equations by taking advantage of the results already existing in the literature. It is in this sense that, inspired by the works  \cite{BT}, \cite{CT12} or \cite{DaSa1} and exploiting these, we present a stochastic wave equation with respect to the Hermite sheet in spatial dimension $d\geq 1$ and we  study the existence, regularity, and other properties of the solution, including the existence of local times and of the joint density.\\

We organize our paper as follows. Section 2 presents the necessary notations and prove several properties of the Hermite sheet. In Section 3, we construct
Wiener integrals with respect to this process. Section 4 is devoted to present the wave equation and discuss the existence and regularity of the solution and other properties.
\\

\section{Notation and the Hermite sheet}

\hspace*{1.1cm} Throughout the work we use the notation introduced in \cite{Breton11}. Fix $d\in \mathbb{N} \backslash \left\lbrace0\right\rbrace$ and consider multi-parametric processes indexed in $\mathbb{R} ^{d}$. We shall use bold notation for multi-indexed quantities, i.e., $\mathbf{a}=(a_{1},a_{2},\ldots ,a_{d})$, $\mathbf{ab}=(a_{1}b_{1},a_{2}b_{2},\ldots ,a_{d}b_{d})$, $\mathbf{a/b}=(a_{1}/b_{1},a_{2}/b_{2},\ldots ,a_{d}/b_{d})$, $ [\mathbf{a},\mathbf{b}]= \displaystyle \prod_{i=1}^{d}[a_{i},b_{i}]$, $(\mathbf{a},\mathbf{b})=\displaystyle \prod_{i=1}^{d}(a_{i},b_{i})$, $\displaystyle\sum_{\mathbf{i}\in [\mathbf{0},\mathbf{N}]} a_{\mathbf{i}} = \displaystyle\sum_{i_{1}=0}^{N_{1}} \displaystyle\sum_{i_{2}=0}^{N_{2}} \ldots \displaystyle\sum_{i_{d}=0}^{N_{d}} a_{i_{1},i_{2},\ldots ,i_{d}}$, $\mathbf{a}^{\mathbf{b}}=\displaystyle \prod_{i=1}^{d} a_{i}^{b_{i}}$, and $ \mathbf{a} < \mathbf{b} $ iff $ a_{1} < b_{1}, a_{2} < b_{2},\ldots ,a_{d} < b_{d}$ (analogously for the other inequalities).\\

Before introducing the {\em Hermite sheet} we briefly recall the {\em fractional Brownian sheet} and the {\em standard Brownian sheet}.\\

The {\em d-parametric anisotropic fractional Brownian sheet} is the centered Gaussian process  $\{B_{\mathbf{t}}^{\mathbf{H}}:\mathbf{t}=(t_{1},\ldots ,t_{d}) \in \mathbb{R} ^{d}\}$ with Hurst multi-index $\mathbf{H}=(H_{1},\ldots ,H_{d}) \in (0,1)^{d}$. It is equal to zero on the hyperplanes $\{\mathbf{t}:t_{i}=0\}$, $1\leq i \leq d$, and its covariance function is given by
\begin{eqnarray}\label{cov-fra}
\nonumber
R_{\mathbf{H}}(\mathbf{s},\mathbf{t}) &=& \mathbb{E}[B_{\mathbf{s}}^{\mathbf{H}}B_{\mathbf{t}}^{\mathbf{H}}] \\
&=& \displaystyle \prod_{i=1}^{d} R_{H_{i}}(s_{i},t_{i}) = \prod_{i=1}^{d} \frac{s_{i}^{2H_{i}} + t_{i}^{2H_{i}} - |t_{i} - s_{i}|^{2H_{i}} }{2}.
\end{eqnarray}
The {\em d-parametric standard Brownian sheet} is the Gaussian process $\{W_{\mathbf{t}}:\mathbf{t}=(t_{1},\ldots ,t_{d}) \in \mathbb{R} ^{d}\}$ equal to zero on the hyperplanes $\{\mathbf{t}:t_{i}=0\}$, $1\leq i \leq d$, and its covariance function is given by
\begin{eqnarray}\label{cov-brow}
R(\mathbf{s},\mathbf{t}) = \mathbb{E}[W_{\mathbf{s}}W_{\mathbf{t}}] = \displaystyle \prod_{i=1}^{d} R(s_{i},t_{i}) = \prod_{i=1}^{d} s_{i} \wedge t_{i}.
\end{eqnarray}

Let $q\in \mathbb{N}$, and the Hurst multi-index $\mathbf{H}=(H_{1},H_{2}, \ldots ,H_{d}) \in (\frac{1}{2},1)^{d}$. The {\em Hermite sheet of order q} is given by

\begin{eqnarray}
\label{hermite-sheet-1}
\nonumber
 Z^{q}_{\mathbf{H}}(\mathbf{t}) &=& c(\mathbf{H},q) \int_{\mathbb{R}^{d\cdot q}}
 \int_{0}^{t_{1}}
 \ldots \int_{0}^{t_{d}} \left( \prod _{j=1}^{q} (s_{1}-y_{1,j})_{+} ^{-\left( \frac{1}{2} + \frac{1-H_{1}}{q} \right) }
 \ldots (s_{d}-y_{d,j})_{+} ^{-\left( \frac{1}{2} + \frac{1-H_{d}}{q} \right) }   \right)  \\
\nonumber
& &  ds_{d} \ldots ds_{1} \quad
dW(y_{1,1},\ldots ,y_{d,1} )\ldots dW(y_{1,q},\ldots ,y_{d,q}) \\
&=&
c(\mathbf{H},q) \int_{\mathbb{R}^{d\cdot q}}
\int_{0}^{\mathbf{t}} \quad \prod _{j=1}^{q} (\mathbf{s}-\mathbf{y}_{j})_{+} ^{-\mathbf{\left( \frac{1}{2} + \frac{1-H}{q} \right)} }  d\mathbf{s} \quad dW(\mathbf{y}_{1})\ldots dW(\mathbf{y}_{q})
\end{eqnarray}
where $x_{+}=\max(x,0)$. For a better understanding about multiple stochastic integrals we refer to \cite{N}. As pointed out before, when $q=1$, (\ref{hermite-sheet-1}) is the fractional Brownian sheet with Hurst multi-index $\mathbf{H}=(H_{1},H_{2}, \ldots ,H_{d}) \in (\frac{1}{2},1)^{d}$. For $q \geq 2$ the process $Z^{q}_{\mathbf{H}}(\mathbf{t})$ is not Gaussian and for $q=2$ we denominate it as the {\em Rosenblatt sheet}.

Now let us calculate the covariance $R_{\mathbf{H}}^{q}(\mathbf{s},\mathbf{t})$ of the Hermite sheet. Using the isometry of multiple Wiener-It\^o integrals and Fubini theorem one gets

\begin{eqnarray*} 
\nonumber
R_{\mathbf{H}}^{q}(\mathbf{s},\mathbf{t}) &=& \mathbb{E}[Z^{q}_{\mathbf{H}}(\mathbf{s})Z^{q}_{\mathbf{H}}(\mathbf{t})] \\
\nonumber
&=& \mathbb{E} \left\lbrace c(\mathbf{H},q)^{2} \int_{\mathbb{R}^{d\cdot q}}
\int_{0}^{\mathbf{s}} \quad \prod _{j=1}^{q} (\mathbf{u}-\mathbf{y}_{j})_{+} ^{-\mathbf{\left( \frac{1}{2} + \frac{1-H}{q} \right)} }  d\mathbf{u} \quad dW(\mathbf{y}_{1})\ldots dW(\mathbf{y}_{q}) \right. \\
\nonumber
& &  \cdot \left. \int_{\mathbb{R}^{d\cdot q}}
\int_{0}^{\mathbf{t}} \quad \prod _{j=1}^{q} (\mathbf{v}-\mathbf{y}_{j})_{+} ^{-\mathbf{\left( \frac{1}{2} + \frac{1-H}{q} \right)} }  d\mathbf{v} \quad dW(\mathbf{y}_{1})\ldots dW(\mathbf{y}_{q}) \right\rbrace \\
\nonumber
&=& c(\mathbf{H},q)^{2} \int_{\mathbb{R}^{d\cdot q}} \left\lbrace \int_{0}^{s_{1}} \ldots \int_{0}^{s_{d}} \prod _{j=1}^{q} \prod _{i=1}^{d} (u_{i}-y_{i,j})_{+} ^{-\left( \frac{1}{2} + \frac{1-H_{i}}{q} \right) }  du_{d} \ldots du_{1} \right. \\
\nonumber
& & \cdot \left. \int_{0}^{t_{1}} \ldots \int_{0}^{t_{d}} \prod _{j=1}^{q} \prod _{i=1}^{d} (v_{i}-y_{i,j})_{+} ^{-\left( \frac{1}{2} + \frac{1-H_{i}}{q} \right) }  dv_{d} \ldots dv_{1} \right\rbrace dy_{1,1}\ldots dy_{d,1}\ldots dy_{1,q} \ldots dy_{d,q} \\
\nonumber
&=& c(\mathbf{H},q)^{2} \int_{0}^{t_{1}} \int_{0}^{s_{1}} \int_{\mathbb{R}^{q}} \prod _{j=1}^{q} (u_{1}-y_{1,j})_{+} ^{-\left( \frac{1}{2} + \frac{1-H_{1}}{q} \right) } (v_{1}-y_{1,j})_{+} ^{-\left( \frac{1}{2} + \frac{1-H_{1}}{q} \right) } dy_{1,1} \ldots dy_{1,q} du_{1} dv_{1} \\
\nonumber
& &  \vdots \\
\nonumber
& & \int_{0}^{t_{d}} \int_{0}^{s_{d}} \int_{\mathbb{R}^{q}} \prod _{j=1}^{q} (u_{d}-y_{d,j})_{+} ^{-\left( \frac{1}{2} + \frac{1-H_{d}}{q} \right) } (v_{d}-y_{d,j})_{+} ^{-\left( \frac{1}{2} + \frac{1-H_{d}}{q} \right) } dy_{d,1} \ldots dy_{d,q} du_{d} dv_{d} \\
\end{eqnarray*}
but
\begin{eqnarray}\label{prop1}
\nonumber
&&\int_{\mathbb{R}^{q}} \prod _{j=1}^{q}  (u-x_{j})_{+} ^{-\left( \frac{1}{2} + \frac{1-H}{q} \right) } (v-x_{j})_{+} ^{-\left( \frac{1}{2} + \frac{1-H}{q} \right) } dx_{1} \ldots dx_{q}   \\
&=&
\left[ \int_{\mathbb{R}} (u-x)_{+} ^{-\left( \frac{1}{2} + \frac{1-H}{q} \right) } (v-x)_{+} ^{-\left( \frac{1}{2} + \frac{1-H}{q} \right) } dx \right] ^{q},
\end{eqnarray}
so
\begin{eqnarray*}
R_{\mathbf{H}}^{q}(\mathbf{s},\mathbf{t}) &=& c(\mathbf{H},q)^{2} \int_{0}^{t_{1}} \int_{0}^{s_{1}} \left[ \int_{\mathbb{R}} (u_{1}-y_{1})_{+} ^{-\left( \frac{1}{2} + \frac{1-H_{1}}{q} \right) } (v_{1}-y_{1})_{+} ^{-\left( \frac{1}{2} + \frac{1-H_{1}}{q} \right) } dy_{1} \right] ^{q} du_{1} dv_{1} \\
& & \vdots \\
\nonumber
& &  \int_{0}^{t_{d}} \int_{0}^{s_{d}} \left[ \int_{\mathbb{R}} (u_{d}-y_{d})_{+} ^{-\left( \frac{1}{2} + \frac{1-H_{d}}{q} \right) } (v_{d}-y_{d})_{+} ^{-\left( \frac{1}{2} + \frac{1-H_{d}}{q} \right) } dy_{d} \right] ^{q} du_{d} dv_{d}. \\
\end{eqnarray*}
Recalling that the Beta function $\beta (p,q)= \int_{0}^{1} z^{p-1}(1-z)^{q-1} dz, p,q>0$, satisfies the following identity
\begin{eqnarray}\label{beta}
 \int_{\mathbb{R}} (u-y)_{+} ^{a-1} (v-y)_{+}^{a-1}
dy& =&\beta (a,1-2a) \vert u-v\vert ^{2a-1}
\end{eqnarray}
if $a\in (0, \frac{1}{2})$, we see that, since $H_{k}\geq \frac{1}{2} $ for every $k=1,..,d$,
\begin{eqnarray*}
R_{\mathbf{H}}^{q}(\mathbf{s},\mathbf{t}) &=& c(\mathbf{H},q)^{2} \int_{0}^{t_{1}} \int_{0}^{s_{1}} \beta \left( \frac{1}{2} - \frac{1-H_{1}}{q}, \frac{2(1-H_{1})}{q} \right) ^{q} \cdot \vert u_{1}-v_{1}\vert ^{2(H_{1}-1)} du_{1} dv_{1} \\
\nonumber
& & \ldots \int_{0}^{t_{d}} \int_{0}^{s_{d}} \beta \left( \frac{1}{2} - \frac{1-H_{d}}{q}, \frac{2(1-H_{d})}{q} \right) ^{q} \cdot \vert u_{d}-v_{d}\vert ^{2(H_{d}-1)} du_{d} dv_{d} \\
\nonumber
&=& c(\mathbf{H},q)^{2} \beta \left( \frac{1}{2} - \frac{1-H_{1}}{q},\frac{2(1-H_{1})}{q} \right) ^{q} \frac{1}{2H_{1}(2H_{1}-1)} \left( s_{1}^{2H_{1}} + t_{1}^{2H_{1}} - \vert t_{1} - s_{1} \vert ^{2H_{1}}   \right) \\
\nonumber
& & \ldots \beta \left( \frac{1}{2} - \frac{1-H_{d}}{q},\frac{2(1-H_{d})}{q} \right) ^{q} \frac{1}{2H_{d}(1-2H_{d})} \left( s_{d}^{2H_{d}} + t_{d}^{2H_{d}} - \vert t_{d} - s_{d} \vert ^{2H_{d}}   \right).
\end{eqnarray*}
So now we choose
\begin{equation}\label{constant}
c(\mathbf{H},q)^{2} = \left( \frac{\beta \left( \frac{1}{2} - \frac{1-H_{1}}{q},\frac{2(1-H_{1})}{q} \right) ^{q}}{H_{1}(2H_{1}-1)} \right)^{-1} \ldots \left( \frac{\beta \left( \frac{1}{2} - \frac{1-H_{d}}{q},\frac{2(1-H_{d})}{q} \right) ^{q}}{H_{d}(2H_{d}-1)} \right)^{-1}
\end{equation}
in this way we get $\mathbb{E} \left( Z^{q}_{\mathbf{H}}(\mathbf{t}) ^{2} \right) = \mathbf{t}^{2\mathbf{H}} = t_{1}^{2H_{1}} \ldots t_{d}^{2H_{d}} $, and finally
\begin{eqnarray}\label{cov-hermite}
\nonumber
R_{\mathbf{H}}^{q}(\mathbf{s},\mathbf{t}) &=& \frac{1}{2}\left( s_{1}^{2H_{1}} + t_{1}^{2H_{1}} - \vert t_{1} - s_{1} \vert ^{2H_{1}}   \right) \ldots \left( s_{d}^{2H_{d}} + t_{d}^{2H_{d}} - \vert t_{d} - s_{d} \vert ^{2H_{d}}   \right) \\
\nonumber
&=& \displaystyle \prod_{i=1}^{d} \frac{s_{i}^{2H_{i}} + t_{i}^{2H_{i}} - |t_{i} - s_{i}|^{2H_{i}} }{2} \\
&=& \prod_{i=1}^{d} R_{H_{i}}(s_{i},t_{i}) = R_{\mathbf{H}}(\mathbf{s},\mathbf{t}).
\end{eqnarray}

\begin{remark}
From the previous development we see that:
\begin{itemize}
\item The covariance structure is the same for all $q\geq 1$, so   it coincides with the covariance of the fractional Brownian sheet.
\item In order to all the quantities are well defined, the condition $H_{k} \in (\frac{1}{2},1)$, $k=1, \ldots , d$ must be satisfied.
\end{itemize} 

\end{remark}

We will next prove the basic properties of the Hermite sheet: self-similarity, stationarity of the increments and H\"older continuity.

Let us first recall the concept of self-similarity for multiparameter stochastic processes.

\begin{definition}\label{self}
A stochastic process $(X_{\mathbf{t} } )_{\mathbf{t}\in T}$, where $T\subset \mathbb{R} ^{d}$  is called self-similar with self-similarity order $\mathbf{\alpha}=(\alpha _{1},\ldots , \alpha _{d})>0$ if for any $\mathbf{h}=(h_{1},\ldots , h_{d})>0$ the stochastic process $(\hat{X} _{\mathbf{t}}) _{\mathbf{t} \in T}$ given by
\begin{equation*}
\hat{X} _{\mathbf{t}} =\mathbf{h} ^{\mathbf{\alpha}} X _{\frac{\mathbf{t}}{\mathbf{h}}} = h_{1} ^{\alpha_{1}} ...h_{d} ^{\alpha _{d}} X_{\frac{t_{1}}{h_{1}},\ldots , \frac{t_{d}}{h_{d}}}
\end{equation*}
has the same finite dimensional distributions as the process $X$.
\end{definition}

\begin{prop}
The Hermite sheet is self-similar of order $\mathbf{H}=(H_{1},\ldots , H_{d}) $.

\end{prop}
{\bf Proof: } The scaling property of the Wiener sheet implies that for every $0<\mathbf{c}= (c_{1},\ldots , c_{d} ) \in \mathbb{R} ^{d}$ the processes $\left( W(\mathbf{c} \mathbf{t} )  _{\mathbf{t} \geq 0} \right) $  and $\left(\sqrt{\mathbf{c}} W( \mathbf{t} ) \right) _{\mathbf{t} \geq 0} $ have the same finite dimensional distributions.
Therefore, if $\mathbf{1}=(1,\ldots ,1) \in \mathbb{R} ^{d}$, using obvious changes of variables in the integrals with respect to $d\mathbf{s}$ and $dW$,
\begin{eqnarray*}
\hat{Z}^{q}_{\mathbf{H}}(t) &=&\mathbf{h} ^{\mathbf{H}}{Z} ^{q} _{\frac{\mathbf{t}}{\mathbf{h}}}\\
&=&
c(\mathbf{H},q)\mathbf{h} ^{\mathbf{H}} \int_{\mathbb{R}^{d\cdot q}}
\int_{0}^{\frac{\mathbf{t}}{\mathbf{h}}} \quad \prod _{j=1}^{q} (\mathbf{s}-\mathbf{y}_{j})_{+} ^{-\mathbf{\left( \frac{1}{2} + \frac{1-H}{q} \right)} }  d\mathbf{s} \quad dW(\mathbf{y}_{1})\ldots dW(\mathbf{y}_{q})\\
&=&c(\mathbf{H},q)\mathbf{h} ^{\mathbf{H}-\mathbf{1}} \int_{\mathbb{R}^{d\cdot q}}
\int_{0}^{\mathbf{t}} \quad \prod _{j=1}^{q} (\frac{\mathbf{s}}{\mathbf{h}} -\mathbf{y}_{j})_{+} ^{-\mathbf{\left( \frac{1}{2} + \frac{1-H}{q} \right)} }  d\mathbf{s} \quad dW(\mathbf{y}_{1})\ldots dW(\mathbf{y}_{q})\\
&=&c(\mathbf{H},q)\mathbf{h} ^{\mathbf{H}-\mathbf{1}} \int_{\mathbb{R}^{d\cdot q}}
\int_{0}^{\mathbf{t}} \quad \prod _{j=1}^{q} (\frac{\mathbf{s}}{\mathbf{h}} - \frac{\mathbf{y}_{j}}{\mathbf{h}})_{+} ^{-\mathbf{\left( \frac{1}{2} + \frac{1-H}{q} \right)} }  d\mathbf{s} \quad dW(\mathbf{h} ^{-1}\mathbf{y}_{1})\ldots dW(\mathbf{h} ^{-1}\mathbf{y}_{q})\\
&=&c(\mathbf{H},q)\mathbf{h} ^{\mathbf{H}-\mathbf{1}}\mathbf{h} ^{q\mathbf{\left( \frac{1}{2} + \frac{1-H}{q} \right)} }
 \int_{\mathbb{R}^{d\cdot q}}
\int_{0}^{\mathbf{t}} \quad \prod _{j=1}^{q} (\mathbf{s}-\mathbf{y}_{j})_{+} ^{-\mathbf{\left( \frac{1}{2} + \frac{1-H}{q} \right)} }  d\mathbf{s} \quad dW(\mathbf{h} ^{-1}\mathbf{y}_{1})\ldots dW(\mathbf{h} ^{-1}\mathbf{y}_{q})\\
&=^{(d)}& c(\mathbf{H},q)\mathbf{h} ^{\mathbf{H}-\mathbf{1}}\mathbf{h} ^{q\mathbf{\left( \frac{1}{2} + \frac{1-H}{q} \right)} }\mathbf{h} ^{-\frac{q}{2}}
 \int_{\mathbb{R}^{d\cdot q}}
\int_{0}^{\mathbf{t}} \quad \prod _{j=1}^{q} (\mathbf{s}-\mathbf{y}_{j})_{+} ^{-\mathbf{\left( \frac{1}{2} + \frac{1-H}{q} \right)} }  d\mathbf{s} \quad dW(\mathbf{y}_{1})\ldots dW(\mathbf{y}_{q})\\
&=&Z^{q}_{\mathbf{H}}(t)
\end{eqnarray*}
where $=^{(d)}$ means equivalence of finite dimensional distributions.

\qed
\\

Let us recall the notion of the increment of a $d$-parameter process $X$ on a rectangle $[\mathbf{s}, \mathbf{t} ] \subset \mathbb{R} ^{d}$, $\mathbf{s}= (s_{1},\ldots , s_{d} ), \mathbf{t}=(t_{1},\ldots ,t_{d})$, with $\mathbf{s} \leq \mathbf{t}$. This increment is denoted by $\Delta X_{[\mathbf{s}, \mathbf{t} ]}$ and it is given by
\begin{equation}\label{mi}
\Delta X_{[\mathbf{s}, \mathbf{t} ]}= \displaystyle\sum_{r\in \{0,1\} ^{d}}(-1) ^{d- \sum_{i=1}^{d} r_{i}} X_{\mathbf{s} + \mathbf{r} \cdot (\mathbf{t}-\mathbf{s})}.
\end{equation}
When $d=1$ one obtains $\Delta X_{[\mathbf{s}, \mathbf{t} ]}=X_{t}-X_{s}$ while for $d=2$ one gets $\Delta X_{[\mathbf{s}, \mathbf{t} ]}=X_{t_{1}, t_{2}} -X_{t_{1}, s_{2}} -X_{s_{1}, t_{2}} + X_{s_{1}, s_{2}}$.

\begin{definition}\label{stationary}
A process $(X_{\mathbf{t}}, \mathbf{t} \in \mathbb{R} ^{d})$ has stationary increments if for every $\mathbf{h} >0, \mathbf{h}\in \mathbb{R} ^{d}$ the stochastic processes $(\Delta X _{[0, \mathbf{t}]}, \mathbf{t} \in \mathbb{R} ^{d})$ and $(\Delta X_{[\mathbf{h}, \mathbf{h}+ \mathbf{t}] },\mathbf{t} \in \mathbb{R} ^{d})$ have the same finite dimensional distributions.

\end{definition}

\begin{prop}
The Hermite sheet $(Z^{q} (\mathbf{t})) _{\mathbf{t} \geq 0}$ has stationary increments.

\end{prop}
{\bf Proof: } Developing the increments of the process using the definition of the Hermite sheet and proceding as in the proof of Proposition 1 using the change of variables $\mathbf{s}' = \mathbf{s} - \mathbf{h}$, it is immediate to see that for every $\mathbf{h} >0, \mathbf{h}\in \mathbb{R} ^{d}$,
$$\Delta Z^{q}_{[\mathbf{h}, \mathbf{h}+ \mathbf{t}] } =^{(d)}  \Delta Z^{q} _{[0, \mathbf{t}]} $$
for every $\mathbf{t}$.
\qed

\begin{prop}
The trajectories of the Hermite sheet $(Z^{q}(\mathbf{t}), \mathbf{t} \geq 0)$ are H\"older continuous of any order $\mathbf{\delta} = (\delta _{1},\ldots , \delta _{d}) \in [0, \mathbf{H} )$ in the following sense: for every $\omega \in \Omega$, there exists a constant $C_{\omega } >0$ such that for every $\mathbf{s}, \mathbf{t} \in \mathbb{R} ^{d}, \mathbf{s}, \mathbf{t}\geq 0$,
\begin{equation*}
\vert \Delta Z^{q} _{[\mathbf{s}, \mathbf{t}]} \vert \leq C_{\omega } \vert t_{1}-s_{1}\vert ^{\delta _{1}}  \cdots \vert t_{d}-s_{d}\vert ^{\delta _{d}}= C _{\omega} \vert \mathbf{t} -\mathbf{s}\vert ^{\mathbf{\delta}}.
\end{equation*}

\end{prop}
{\bf Proof: } Using the Cencov's criteria (see \cite{Cencov}) and the fact that the process $Z^{q}$ is almost surely equal to 0 when $t_{i}=0$, it suffices to check that
\begin{equation}
\label{i1}
\mathbb{E} \left| \Delta Z^{q} _{[\mathbf{s}, \mathbf{t}]} \right| ^{p} \leq C\left( \vert t_{1}-s_{1}\vert \cdots  \vert t_{d}-s_{d}\vert \right) ^{1+\gamma }
\end{equation}
for some $p\geq 2$ and $\gamma >0$. From the self-similarity and the stationarity of the increments of the process $Z^{q}$, we have for every $p\geq 2$
\begin{equation*}
\mathbb{E} \left| \Delta Z^{q} _{[\mathbf{s}, \mathbf{t}]} \right| ^{p} =\mathbb{E} \left| Z_{\mathbf{1}} \right| ^{p} \left( \vert t_{1}-s_{1}\vert \cdots  \vert t_{d}- s_{d}\vert \right)^{p\mathbf{H}}
\end{equation*}
and this obviously implies (\ref{i1}) by taking $p$ arbitrary large.
\qed
\\

\begin{remark}
In the one-parameter case, there exists several representations of the Hermite process (spectral domain representation, finite interval representation, positive half-axis representation, time domain representation, see \cite{PT10-Hp}). It has been shown in \cite{PT10-Hp} that all these representations of the (one-parameter) Hermite process have the same finite dimensional distributions. It would be interesting to generalize this study to the multi-parameter case.

\end{remark}

\section{Wiener integrals with respect to the Hermite sheet}

\hspace*{1.1cm} Now we are well positioned to present Wiener integrals with respect to the $d$-parametric Hermite sheet. Let us consider a Hermite sheet $\left( Z^{q} _{\mathbf{H}} (\mathbf{t})\right) _{\mathbf{t}\in \mathbb{R} ^{d}}.$  Denote $\mathscr{E}$ the family of elementary functions on $\mathbb{R}^{d}$ of the form

\begin{eqnarray}\label{elem-func}
 f(\mathbf{u}) &=&\sum _{l=1}^{n} a_{l}1_{(\mathbf{t}_{l}, \mathbf{t}_{l+1}] }(\mathbf{u}) \\
 \nonumber
&=& \sum _{l=1}^{n} a_{l} 1_{({t}_{1,l}, {t}_{1,l+1}] \times \ldots \times ({t}_{d,l}, {t}_{d,l+1}]} (u_{1}, \ldots ,u_{d}),
\hskip0.5cm  \mathbf{t}_{l} < \mathbf{t}_{l+1}, \hskip0.2cm a_{l} \in \mathbb{R},\hskip0.2cm l=1,\ldots ,n.
\end{eqnarray}
For functions like $f$ above we can naturally define its Wiener integral with respect to the Hermite sheet $Z^{q}_{\mathbf{H}}$ as
\begin{eqnarray} \label{Her-int-1}
 \int_{\mathbb{R}^{d}}  f(\mathbf{u}) dZ_{\mathbf{H}}^{q} (\mathbf{u})&=& \displaystyle\sum_{l=1}^{n} a_{l} (\Delta Z^{q} _{\mathbf{H}}) _{ [\mathbf{t} _{l}, \mathbf{t} _{l+1} ]}
 \end{eqnarray}
where $(\Delta  Z^{q} _{\mathbf{H}}) _{ [\mathbf{t} _{l}, \mathbf{t} _{l+1} ]}$ (see (\ref{mi})) stands for the generalized increments of $Z^{q}_{\mathbf{H}}$ on the rectangle
\begin{equation*}
\Delta_{\mathbf{t}_{l}} := \left[ \mathbf{t}_{l}, \mathbf{t}_{l+1} \right] = \displaystyle \prod_{i=1}^{d} \left[t_{i,l}, t_{i,l+1} \right].
\end{equation*}
In the case $d=1$, we simply have
$$(\Delta  Z^{q} _{\mathbf{H}}) _{[ \mathbf{t} _{l}, \mathbf{t} _{l+1} ]}= Z^{q} _{\mathbf{H} } (t_{1,l+1} -t_{1,l}) $$
while for $d=2$
$$(\Delta  Z^{q} _{\mathbf{H}}) _{[ \mathbf{t} _{l}, \mathbf{t} _{l+1} ]}= Z^{q}_{\mathbf{H}} (t_{1,l+1}, t_{2, l+1} ) -  Z^{q}_{\mathbf{H}} (t_{1,l}, t_{2, l+1} )- Z^{q}_{\mathbf{H}} (t_{1,l+1}, t_{2, l} )+ Z^{q}_{\mathbf{H}} (t_{1,l}, t_{2, l} ).$$

With the purpose of extending the definition (\ref{Her-int-1}) to a larger family of integrands, we will point out some observations before. Let us consider the mapping $J$ on the set of functions $f:\mathbf{R}^{d} \rightarrow \mathbf{R}$ to the set of functions $f:\mathbf{R}^{d\cdot q} \rightarrow \mathbf{R}$ such that
\begin{eqnarray}\label{J}
J(f)(\mathbf{y}_{1}, \ldots ,\mathbf{y}_{q}) &=&  c(\mathbf{H},q) \int_{\mathbb{R}^{d}} f(\mathbf{u})  \prod_{j=1}^{q} (\mathbf{u} - \mathbf{y}_{j})_{+}^{-\mathbf{ \left( \frac{1}{2} + \frac{1-H}{q} \right)}} d\mathbf{u} \\
\nonumber
&=&
c(\mathbf{H},q) \int_{\mathbb{R}^{d}} f(u_{1}, \ldots ,u_{d})  \prod_{j=1}^{q} \displaystyle \prod_{i=1}^{d} (u_{i}-u_{i,j})_{+}^{-\left( \frac{1}{2} + \frac{1-H_{i}}{q} \right)} du_{1} \ldots du_{d}.
\end{eqnarray}
Using the mapping $J$ we see that definition (\ref{hermite-sheet-1}) can be re-expressed as follows
\begin{eqnarray}\label{hermite-sheet-2}
 Z^{q}_{\mathbf{H}}(\mathbf{t}) 
&=& \int_{\mathbb{R}^{d\cdot q}}
J(1_{[0,t_{1}]\times \ldots \times[0,t_{d}]}) (\mathbf{y}_{1}, \ldots ,\mathbf{y}_{q})  dW(\mathbf{y}_{1})\ldots dW(\mathbf{y}_{q}).
\end{eqnarray}
Since $J$ is clearly linear, definition (\ref{Her-int-1}) can be tailored to
\begin{eqnarray}\label{Her-int-2}
\nonumber
\int_{\mathbb{R}^{d}}  f(\mathbf{u}) dZ_{\mathbf{H}}^{q} (\mathbf{u})&=& \displaystyle\sum_{l=1}^{n} a_{l} \left( \Delta Z^{q}_{\mathbf{H}}\right)_ {[\mathbf{t}_{l},\mathbf{t}_{l+1}] } \\
\nonumber
&=& 
\displaystyle\sum_{l=1}^{n} a_{l} \left( \displaystyle\sum_{\mathbf{\xi} \in \{0,1\}^{d}} (-1)^{d- \sum_{i=1}^{d} \xi_{i} } Z^{q}_{\mathbf{H}}(t_{1,l+\xi_{1}}, \ldots , t_{d,l+\xi_{d}}) \right) \\
\nonumber
&=& 
\displaystyle\sum_{l=1}^{n} a_{l} \displaystyle\sum_{\mathbf{\xi} \in \{0,1\}^{d}} (-1)^{d- \sum_{i=1}^{d} \xi_{i} } \int_{\mathbb{R}^{d\cdot q}} J(1_{[0,t_{1,l+\xi_{1}}]\times \ldots \times[0,t_{d,l+\xi_{d}]}}) (\mathbf{y}_{1}, \ldots ,\mathbf{y}_{q})  dW(\mathbf{y}_{1})\ldots  dW(\mathbf{y}_{q}) \\
\nonumber
&=& \displaystyle\sum_{l=1}^{n} a_{l} \int_{\mathbb{R}^{d\cdot q}} J(1_{[t_{1,l},t_{1,l+1}]\times \ldots \times[t_{d,l},t_{d,l+1}]}) (\mathbf{y}_{1}, \ldots ,\mathbf{y}_{q})  dW(\mathbf{y}_{1})\ldots  dW(\mathbf{y}_{q}) \\
&=& \int_{\mathbb{R}^{d\cdot q}} J(f) (\mathbf{y}_{1}, \ldots ,\mathbf{y}_{q})  dW(\mathbf{y}_{1})\ldots  dW(\mathbf{y}_{q}).
\end{eqnarray}
In this way we introduce the space
\begin{equation}\label{Integrands-space-1}
\mathcal{H} = \left\lbrace f:\mathbb{R}^{d} \rightarrow \mathbb{R} : \int_{\mathbb{R}^{d\cdot q}} \left( J(f) (\mathbf{y}_{1}, \ldots ,\mathbf{y}_{q}) \right)^{2} d\mathbf{y}_{1} \ldots d\mathbf{y}_{q} < \infty \right\rbrace
\end{equation}
equipped with the norm
\begin{equation}\label{norm-integrands-space-1}
\Vert f\Vert _{\mathcal{H}}^{2}  = \int_{\mathbb{R}^{d\cdot q}}\left( J(f) (\mathbf{y}_{1}, \ldots ,\mathbf{y}_{q}) \right)^{2} d\mathbf{y}_{1} \ldots d\mathbf{y}_{q}.
\end{equation}
Working the expression for the norm we see that
\begin{eqnarray}
\nonumber
\Vert f\Vert _{\mathcal{H}}^{2} &=& c(\mathbf{H},q)^{2} \int_{\mathbb{R}^{d\cdot q}} \left\lbrace \left(  \int_{\mathbb{R}^{d}} f(\mathbf{u}) \prod_{j=1}^{q} (\mathbf{u} - \mathbf{y}_{j})_{+}^{-\mathbf{ \left( \frac{1}{2} + \frac{1-H}{q} \right)}} d\mathbf{u} \right) \right. \\
\nonumber
& & \cdot \left. \left(  \int_{\mathbb{R}^{d}} f(\mathbf{v}) \prod_{j=1}^{q} (\mathbf{v} - \mathbf{y}_{j})_{+}^{-\mathbf{ \left( \frac{1}{2} + \frac{1-H}{q} \right)}} d\mathbf{v} \right) \right\rbrace d\mathbf{y}_{1} \ldots d\mathbf{y}_{q}. 
\end{eqnarray}
Using (\ref{prop1}), (\ref{beta}) and (\ref{constant}) we get that
\begin{eqnarray} \label{norm-integrands-space-2}
\nonumber
\Vert f\Vert _{\mathcal{H}}^{2} &=& c(\mathbf{H},q)^{2} \int_{\mathbb{R}^{d}} \int_{\mathbb{R}^{d}} f(u_{1}, \ldots ,u_{d}) f(v_{1}, \ldots ,v_{d}) \\
\nonumber
& & \left\lbrace \displaystyle \prod_{i=1}^{d} \int_{\mathbb{R}^{q}} \prod_{j=1}^{q} (u_{i}-y_{i,j})_{+}^{-\left( \frac{1}{2} + \frac{1-H_{i}}{q} \right)} (v_{i}-y_{i,j})_{+}^{-\left( \frac{1}{2} + \frac{1-H_{i}}{q} \right)} dy_{i,1} \ldots dy_{i,q}  \right\rbrace du_{1} \ldots du_{d} \ \  dv_{1} \ldots dv_{d} \\
\nonumber
&=&  c(\mathbf{H},q)^{2} \int_{\mathbb{R}^{d}} \int_{\mathbb{R}^{d}} f(u_{1}, \ldots ,u_{d}) f(v_{1}, \ldots ,v_{d}) \\
\nonumber
& & \cdot \displaystyle \prod_{i=1}^{d} \left( \int_{\mathbb{R}} (u_{i}-y)_{+}^{-\left( \frac{1}{2} + \frac{1-H_{i}}{q} \right)} (v_{i}-y)_{+}^{-\left( \frac{1}{2} + \frac{1-H_{i}}{q} \right)} dy \right) ^{q} du_{1} \ldots du_{d} \ \  dv_{1} \ldots dv_{d} \\
\nonumber
&=&  \int_{\mathbb{R}^{d}} \int_{\mathbb{R}^{d}} f(u_{1}, \ldots ,u_{d}) f(v_{1}, \ldots ,v_{d}) \displaystyle \prod_{i=1}^{d} H_{i}(2H_{i}-1) \vert u - v\vert ^{2H_{i}-2} du_{1} \ldots du_{d} \ \ dv_{1} \ldots dv_{d} \\
&=&
\mathbf{H(2H-1)} \int_{\mathbb{R}^{d}} \int_{\mathbb{R}^{d}} f(\mathbf{u}) f( \mathbf{v}) \vert \mathbf{u} - \mathbf{v}\vert ^{\mathbf{2H-2}} d\mathbf{u} d\mathbf{v},
\end{eqnarray}
hence
\begin{equation}\label{Integrands-space-2}
\mathcal{H}=\left\lbrace f:\mathbb{R}^{d} \rightarrow \mathbb{R}: \int_{\mathbb{R}^{d}} \int_{\mathbb{R}^{d}} f(\mathbf{u}) f( \mathbf{v}) \vert \mathbf{u} - \mathbf{v}\vert ^{\mathbf{2H-2}} d\mathbf{u} d\mathbf{v} < + \infty  \right\rbrace
\end{equation}
and
\begin{equation*}
\Vert f\Vert _{\mathcal{H}}^{2} = \mathbf{H(2H-1)}\int_{\mathbb{R}^{d}} \int_{\mathbb{R}^{d}} f(\mathbf{u}) f( \mathbf{v}) \vert \mathbf{u} - \mathbf{v}\vert ^{\mathbf{2H-2}} d\mathbf{u} d\mathbf{v}.
\end{equation*}

The mapping
\begin{equation}
f \rightarrow \int_{\mathbb{R}^{d}} f(\mathbf{u}) dZ^{q}_{\mathbf{H}}(\mathbf{u})
\end{equation}
provides an isometry from $\mathscr{E}$ to $L^{2}(\Omega)$. Indeed, for $f$ like in (\ref{elem-func}) it holds that
\begin{eqnarray}\label{Isometry-1}
\nonumber
&&\mathbb{E}\left\lbrace \left( \int_{\mathbb{R}^{d}} f(\mathbf{u}) dZ^{q}_{\mathbf{H}}(\mathbf{u}) \right)^{2} \right\rbrace \\
&=&
 \displaystyle \sum_{k,l=0}^{n-1} a_{k}a_{l} \mathbb{E}\left( \left( \Delta  Z^{q}_{\mathbf{H}} \right)_{[\mathbf{t}_{k},\mathbf{t}_{k+1}]} \cdot \left( \Delta Z^{q}_{\mathbf{H}}\right)_ {[\mathbf{t}_{l},\mathbf{t}_{l+1}] } \right) \\
\nonumber
&=&
 \displaystyle\sum_{k,l=0}^{n-1} a_{k}a_{l} \displaystyle\sum_{\mathbf{\xi} \in \{0,1\}^{d}} (-1)^{d- \sum_{i=1}^{d} \xi_{i} } \displaystyle\sum_{\mathbf{\rho} \in \{0,1\}^{d}} (-1)^{d- \sum_{j=1}^{d} \rho_{j} } \mathbb{E} \left\lbrace Z^{q}_{\mathbf{H}}(\mathbf{t}_{k+\xi}) Z^{q}_{\mathbf{H}}(\mathbf{t}_{l+\rho}) \right\rbrace \\
\nonumber
&=&
 \displaystyle\sum_{k,l=0}^{n-1} a_{k}a_{l} \displaystyle\sum_{\mathbf{\xi} \in \{0,1\}^{d}} (-1)^{d- \sum_{i=1}^{d} \xi_{i} } \displaystyle\sum_{\mathbf{\rho} \in \{0,1\}^{d}} (-1)^{d- \sum_{j=1}^{d} \rho_{j} }  R_{\mathbf{H}}(\mathbf{t}_{k+\xi},\mathbf{t}_{l+\rho}) \\
\nonumber
&=& \displaystyle\sum_{k,l=0}^{n-1} a_{k}a_{l} H_{1}(2H_{1}-1) \ldots H_{d}(2H_{d}-1) \int_{t_{1,k}}^{t_{1,k+1}} \ldots \int_{t_{d,k}}^{t_{d,k+1}} \cdot \int_{t_{1,l}}^{t_{1,l+1}} \ldots \int_{t_{d,l}}^{t_{d,l+1}} \\
\nonumber
& & \vert u_{1} - v_{1} \vert ^{2H_{1}-2} \ldots \vert u_{d} - v_{d} \vert ^{2H_{d}-2} du_{1} \ldots du_{d} dv_{1} \ldots dv_{d} \\
\nonumber
&=& \displaystyle\sum_{k,l=0}^{n-1} a_{k}a_{l} <1_{[t_{1,k},t_{1,k+1}]\times \cdots \times [t_{d,k},t_{d,k+1}]},1_{[t_{1,l},t_{1,l+1}]\times \cdots \times [t_{d,l},t_{d,l+1}]}>_{\mathcal{H}} \\
&=& <f,f>_{\mathcal{H}},
\end{eqnarray}
where we have made a slight abuse of notation, $\mathbf{t}_{k+\xi} = (t_{1,k+\xi_{1}}, \ldots , t_{d,k+\xi_{d}})$.

\vspace*{0.5cm}

 On the other hand, from what shown in \cite{PiTa1} section 4, it follows that the set of elementary functions
${\cal{E}}$ is dense in  ${\cal{H}}$.  As a consequence the
mapping (\ref{J}) can be extended to an isometry from
${\cal{H}}$ to $L^{2}(\Omega)$ and relation (\ref{hermite-sheet-2}) still
holds.

\vskip0.1cm

\begin{remark}
The elements of ${\cal{H}}$ may be not functions
but distributions; it is therefore more practical to work with
subspaces of ${\cal{H}}$ that are sets of functions. Such a subspace
is
\begin{eqnarray*}  \left| {\cal{H}}\right| &=&\left \{ f:\mathbb{R}^{d}\to
\mathbb{R} \,\,\Big |  \int _{\mathbb{R}^{d}} \int_{\mathbb{R}^{d}} \vert f(\mathbf{u})\vert
\vert f(\mathbf{v})\vert \vert \mathbf{u}-\mathbf{v}\vert ^{\mathbf{2H-2}} d\mathbf{v}d\mathbf{u} <\infty \right \}.
\end{eqnarray*}
Then $\left| {\cal{H}}\right|$ is a strict subspace of $
{\cal{H}}$ and we actually have the inclusions
\begin{eqnarray}
\label{inclu1} L^{\mathbf{2}}(\mathbb{R}^{d}) \cap  L^{\mathbf{1}}(\mathbb{R}^{d})\subset
L^{\frac{1}{\mathbf{H} }} (\mathbb{R}^{d}) \subset \left| {\cal{H}}\right|
\subset  {\cal{H}},
\end{eqnarray}
where $L^{\mathbf{p}}$ denotes $L^{p_{1}} \otimes \ldots \otimes L^{p_{d}} $. \\
 The space $\left| {\cal{H}}\right|$ is not complete with respect to the norm $\Vert \cdot
\Vert _{{\cal{H}}}$ but it is a Banach  space with respect to the
norm
\begin{eqnarray*}
 \Vert f\Vert ^{2}_{\left| {\cal{H}}\right|
}&=& \int _{\mathbb{R}^{d}} \int_{\mathbb{R}^{d}} \vert f(\mathbf{u})\vert
\vert f(\mathbf{v})\vert \vert \mathbf{u}-\mathbf{v}\vert ^{\mathbf{2H-2}} d\mathbf{v}d\mathbf{u}
\end{eqnarray*}

\end{remark}

\begin{remark}

 Expression (\ref{Her-int-2}) presents a useful interpretation for the Wiener integrals with respect to the Hermite sheet; as elements in the $q$-th Wiener chaos generated by the $d$-parametric standard Brownian field.

\end{remark}

\section{Application: The Hermite stochastic wave equation}

\hspace*{1.1cm} In this section we presents the linear stochastic wave equation as an example of equations driven by a Hermite sheet. We show the existence of the solution and study some properties of it thanks to the definition of the Wiener integrals with respect to the Hermite sheet.\\

Consider the linear stochastic wave equation driven by an infinite-dimensional Hermite sheet $ Z^{q}_{\mathbf{H}} $ with Hurst multi-index $\mathbf{H}\in (1/2,1)^{(d+1)}$. That is
\begin{eqnarray}
\label{wave}
 \frac{\partial^2 u}{\partial t^2}(t,\mathbf{x})&=& \Delta u
(t,\mathbf{x}) +\dot Z^{q}_{\mathbf{H}}(t,\mathbf{x}), \quad t>0, \mathbf{x} \in \mathbb{R}^{d}  \\
\nonumber
u(0,\mathbf{x})&=& 0, \quad \mathbf{x} \in \mathbb{R}^{d} \\
\nonumber
\frac{\partial u}{\partial t}(0,\mathbf{x}) &=& 0, \quad \mathbf{x} \in \mathbb{R}^{d}.
\end{eqnarray}
Here $\Delta$ is the Laplacian on $\mathbb{R} ^{d}$ and $Z^{q}_{\mathbf{H}}=\{Z^{q}_{\mathbf{H}}(t,\mathbf{x}); t \geq 0, \mathbf{x} \in \mathbb{R} ^{d}\}$ is the $(d+1)$-parametric Hermite sheet whose covariance is given by
\begin{equation*}
\mathbb{E} \left\lbrace  Z^{q}_{\mathbf{H}}(s,\mathbf{x})  Z^{q}_{\mathbf{H}}(t,\mathbf{y}) \right\rbrace =R_{H}(t,s) R_{\mathbf{H}_{0}} (\mathbf{x}, \mathbf{y} )
\end{equation*}
if $\mathbf{H} = (H, H_{1},\ldots , H_{d})$ and we denoted by $\mathbf{H}_{0}= (H_{1},\ldots , H_{d})$ and $\dot Z^{q} $ stands for the formal derivative of $Z^{(q)}$. Equivalently  we can write
\begin{equation}
\mathbb{E} \left\lbrace \dot Z^{q}_{\mathbf{H}}(s,\mathbf{x}) \dot Z^{q}_{\mathbf{H}}(t,\mathbf{y}) \right\rbrace = H(2H-1)|t-s|^{2H-2}\displaystyle \prod_{i=1}^{d}(H_{i}(2H_{i}-1) \cdot |x_i-y_i|^{2H_i-2}).
\end{equation}

Let $G_1$ be the fundamental solution of $\frac{\partial^2 u}{\partial t^2}-\Delta u=0$. It is known that $G_1(t, \cdot)$ is a distribution in ${\cal{S}'}(\mathbb{R}^d)$ with
rapid decrease, and
\begin{equation}
\label{Fourier-G-wave} \cF
G_1(t,\cdot)(\xi)=\frac{\sin(t|\xi|)}{|\xi|},
\end{equation}
for any $\xi \in \mathbb{R}^d$, $t>0$ and $d \geq 1$, where $\cF$ denotes the Fourier transform (see e.g. \cite{treves75}). In particular,
\begin{eqnarray*}
G_1(t,\mathbf{x})&=&\frac{1}{2}1_{\{|\mathbf{x}|<t\}}, \quad \mbox{if} \ d=1 \\
G_1(t,\mathbf{x})&=&\frac{1}{2 \pi}\frac{1}{\sqrt{t^2-|\mathbf{x}|^2}}1_{\{|\mathbf{x}|<t\}},
\quad \mbox{if} \ d=2 \\
G_1(t,\mathbf{x})&=&c_{d}\frac{1}{t}\sigma_t, \quad \mbox{if} \ d=3,
\end{eqnarray*}
where $\sigma_t$ denotes the surface measure on the 3-dimensional
sphere of radius $t$.

The {\it mild} solution of (\ref{wave}) is a square-integrable process
$u=\{u(t,\mathbf{x}); t \geq 0, \mathbf{x} \in \mathbb{R} ^{d}\}$ defined by:
\begin{equation} \label{def-sol-wave} u(t,\mathbf{x})=\int_{0}^{t}
\int_{\mathbb{R} ^{d}}G_1(t-s,\mathbf{x}-\mathbf{y})Z^{q}_{\mathbf{H}}(ds,d\mathbf{y}).
\end{equation}
The above integral is a Wiener integral with respect to the Hermite sheet, as introduced in Section 2.

\subsection{Existence and regularity of the solution}

By definition, $u(t,\mathbf{x})$ exists if and only if the stochastic integral above is well-defined, i.e. $g_{t,\mathbf{x}}:=G_1(t-\cdot,\mathbf{x}-\cdot)
\in \mathcal{H}$. In this case, $\mathbb{E}|u(t,x)|^2 = \|g_{t,\mathbf{x}}\|_{\cH}^2$.

We  state the result on the existence and the regularity of the solution of (\ref{wave}).
\begin{prop}
Let $Z^{q}_{\mathbf{H}}(t,\mathbf{x})$ be the $(d+1)$-parametric Hermite sheet of order $q$. Denote by
\begin{equation}
\label{beta1}
\beta = d-\displaystyle\sum_{i=1}^{d} (2H_{i}-1).
\end{equation} 
 Then, the stochastic wave equation (\ref{wave}) admits an unique mild solution $(u(t,\mathbf{x}) )_{t\geq 0, \mathbf{x} \in \mathbb{R} ^{d} } $ if and only if
\begin{equation}\label{exist-sol}
    \beta < 2H + 1.  	
\end{equation}
If in addition we have that 
\begin{equation}\label{beta-1}
\beta \in (2H-1, d\wedge 2H+1)
\end{equation}
then, for fixed $0<t_{0}<T$, the following statements are true:
\begin{description}

\item[a.-] Let $\mathbf{x} \in \mathbb{R} ^{d}$ fixed. Then there exist positive constants $c_{1}, c_{2}$ such that for every $s,t\in [t_{0}, T]$
\begin{equation*}
c_{1} \vert t-s\vert ^{2H+1-\beta} \leq \mathbb{E} \left| u(t,\mathbf{x})-u(s,\mathbf{x}) \right| ^{2} \leq c_{2} \vert t-s\vert  ^{2H+1-\beta}.
\end{equation*}
Also for every fixed $\mathbf{x} \in \mathbb{R} ^{d}$  the application
$$t\to u(t,\mathbf{x})$$
is almost surely  H\"older continuous of any order $\delta \in \left( 0, \frac{2H+1-\beta}{2}\right). $ \\

\item[b.-] Fix $M>0$ and $t\in [t_{0}, T]$. Then there exist positive constants $c_{3}, c_{4}$ such that for any $ \mathbf{x}, \mathbf{y} \in [-M,M]^{d} $
\begin{equation*}
c_{3} \vert \mathbf{x} - \mathbf{y} \vert ^{2H+1-\beta} \leq \mathbb{E} \left| u(t,\mathbf{x}) -u(t,\mathbf{y}) \right| ^{2} \leq c_{4} \vert \mathbf{x} - \mathbf{y} \vert ^{2H+1-\beta }.
\end{equation*}
Also, for any $t\in [t_{0}, T]$ the application
$$ \mathbf{x} \to u(t,\mathbf{x})$$
is almost surely H\"older continuous of any order $\delta \in \left(0, \left( \frac{2H+1-\beta}{2} \right) \wedge 1 \right)$.
\item[c.-] Denote by $\Delta $ the following metric on $[0,T] \times \mathbb{R} ^{d}$
\begin{equation}
\label{metric}
\Delta \left( (t,\mathbf{x}); (s,\mathbf{y}) \right)= \vert t-s\vert ^{2H+1-\beta } + \vert \mathbf{x}-\mathbf{y}\vert^{ 2H+1-\beta }.
\end{equation}
Fix $M>0$ and assume (\ref{exist-sol}). For every $t,s \in [t_{0}, T]$ and $\mathbf{x}, \mathbf{y} \in  [-M,M]^{d}$ there exist positive constants $C_{1}, C_{2}$ such that
\begin{equation}\label{increm}
C_{1}\Delta \left( (t,\mathbf{x}); (s,\mathbf{y}) \right) \leq \mathbb{E} \left| u(t,\mathbf{x}) -u(s,\mathbf{y}) \right| ^{2} \leq C_{2} \Delta \left( (t,\mathbf{x}); (s,\mathbf{y}) \right).
\end{equation}
\end{description}

\end{prop}
{\bf Proof: } By the isometry of the Wiener integral with respect to the Hermite sheet, the $L^{2}$ norm  will be
\begin{eqnarray*}
\mathbb{E} u(t, \mathbf{x} ) ^{2}&=& \alpha _{H} \int_{0} ^{t} du \int_{0}^{t} dv \vert u-v \vert ^{2H-2}\int_{\mathbb{R} ^{d}} \int_{\mathbb{R} ^{d}}d\mathbf{y} d\mathbf{z}  G_{1} (t-u, \mathbf{x} -\mathbf{y} ) G_{1} (t-v, \mathbf{x} -\mathbf{z} )\\
 &&\times \displaystyle \prod_{i=1}^{d} (H_{i} (2H_{i}-1)) \vert y_{i} -z_{i} \vert ^{2H_{i} -2}\\
&=& \alpha _{H} \int_{0} ^{t} du \int_{0}^{t} dv \vert u-v \vert ^{2H-2} \int_{\mathbb{R} ^{d}} \frac{\sin (u\vert \mathbf{\xi} \vert ) \sin (v\vert \mathbf{\xi} \vert ) }{\vert \mathbf{\xi} \vert ^{2}} \mu (d\mathbf{\xi} )
\end{eqnarray*}
where
\begin{equation}
\label{mu1}
\mu (d\mathbf{\xi})= c_{\mathbf{H}} \displaystyle \prod_{i=1}^{d}\vert \xi _{i} \vert ^{-(2H_{i} -1)}
\end{equation}
 with $\mathbf{\xi} = (\xi_{1},\ldots , \xi_{d})$, $\alpha _{H}= H(2H-1)$ and $c_{\mathbf{H}}=\mathbf{H}(2\mathbf{H}-1)$. This is, $u(t,\mathbf{x})$ has the same $L^{2}$ norm as in the case $q=1$, i.e. when the noise of the equation is a fractional Brownian sheet. It therefore follows from \cite{BT}, Theorem 3.1 that the above integral is finite if and only if
 \begin{equation*}
 \int_{\mathbb{R} ^{d}} \left( \frac{1}{1+ \vert \mathbf{\xi} \vert ^{2}}\right) ^{H+ \frac{1}{2}} \mu (d\mathbf{\xi} ) < \infty
 \end{equation*}
 with $\mu$ given by (\ref{mu1}). The above condition is equivalent to $\beta < 2H+1$, see Example 3.4 in \cite{BT}.

 The proof of the other items is strongly held in the covariance structure of the Hermite sheet, which is the same as for the fractional Brownian sheet. By a careful revision of the proofs of Theorem 3.1 in \cite{BT}, Propositions 1, 2, 3 and Corollary 1 in \cite{CT12}, it is possible to appreciate that the computations are also valid for any process with a covariance structure like the one presented in these articles, in particular in our case.
\begin{description}
\item{$\bullet $} The bounds for the increments in time are consequence of Proposition 1 in \cite{CT12}, and the H\"older regularity comes from Corollary 1 in \cite{CT12}, this proves \textbf{a}.
\item{$\bullet$ } The bounds for the increments in the space variable are deduced from Proposition 2 in \cite{CT12}, and the space H\"older regularity is direct from Proposition 3 in \cite{CT12}, this proves \textbf{b}.
  \item{$\bullet$ } Point $\mathbf{c}$ follows from \textbf{a} and \textbf{b} by following the lines of  the proof of Theorem 2 in \cite{CT12}.
\end{description}
\qed
\\

\subsection{Existence of local times}

\hspace*{1cm} We will show that the solution of (\ref{wave}), viewed as a process in $(t,\mathbf{x})$, admits a  square integrable local time.

 Let us define the local time of a stochastic process $(X_{t})_{t\in J}$. Here $J$ denotes a subset of $\mathbb{R} ^{d}$.  For any Borel set $I\subset J$ the occupation measure of $X$ on $I$ is defined as
 $$\mu _{I} (A)= \lambda \left( \mathbf{t} \in I, X_{t} \in A\right), \hskip0.5cm A \in {\cal{B}}({\mathbb{R}})$$
 where $ \lambda $ denotes the Lebesgue measure. If $\mu _{I}$ is absolutely continuous with respect to the Lebesgue measure, we say that $X$ has local time on $I$. The local time is defined as the Radon-Nikodym derivative of $\mu _{I}$
 \begin{equation*}
 L(I, x) = \frac{ d\mu _{I}}{d\lambda } (x), \hskip0.5cm x\in \mathbb{R}.
 \end{equation*}
 We will use the notation
\begin{equation*}
L(\mathbf{t},x):= L([0,\mathbf{t}], x), \hskip0.5cm \mathbf{t}\in \mathbb{R}^{d} _{+}, x\in \mathbb{R}.
\end{equation*}
The local time satisfies the occupation time formula
\begin{equation}
\label{oc}
\int_{I} f(X_{\mathbf{t}}) d\mathbf{t} = \int_{\mathbb{R}} f(y) L(I, y) dy
\end{equation}
for any  Borel set $I$ in $T$ and for any measurable function $f: \mathbb{R} \to \mathbb{R}$.

\begin{prop}
Let $u(t,\mathbf{x}), t\geq 0, \mathbf{x} \in \mathbb{R} ^{d}$ be the solution to (\ref{wave}) and assume (\ref{beta-1}) where $\beta $ is given by (\ref{beta1}). Then on each set $[a,b] \times [A,B] \subset [0, \infty) \times \mathbb{R} ^{d} $  the process $\left( u(t,\mathbf{x}),  t\geq 0, \mathbf{x} \in \mathbb{R} ^{d}\right) $ admits a local time $\left( L([a,b] \times [A,B], y), y\in \mathbb{R}\right)$
which is square integrable with respect to $y$
\begin{equation*}
\mathbb{E} \int_{\mathbb{R}}  L([a,b] \times [A,B], y)^{2}dy <\infty \mbox{ a.s. }
\end{equation*}
\end{prop}
{\bf Proof: }It is well known from \cite{B2} (see also Lemma 8.1 in \cite{Xiao1}) that, for a jointly measurable zero-mean stochastic process $X=(X(\mathbf{t}), \mathbf{t}\in[0,\mathbf{T}])$ ($\mathbf{T}$ belongs to $\mathbb{R} ^{d}$) with bounded variance, the condition
$$\int_{[0, \mathbf{T}]}\int_{[0, \mathbf{T}]}(\mathbb{E}[X(\mathbf{t})-X(\mathbf{s})]^{2})^{-1/2}d\mathbf{s} d\mathbf{t}<\infty$$
is sufficient for the local time  of $X$ to exist on $[0,\mathbf{T}]$ almost surely and to be square integrable as a function of $y$. \\
 According to the inequality (\ref{increm}),   for all $I=[a,b]\times [A,B]$ included in $ [0, \infty) \times \mathbb{R} ^{d} $  we have,
$$\int_{I} \int_{I} (\mathbb{E} \left[ u(t,\mathbf{x})-u(s, \mathbf{y}) \right]^{2} )^{-1/2}dtd\mathbf{x} ds d\mathbf{y} < C \int_I \int_I \left( |t-s|^{2H+1-\beta }+ \vert \mathbf{x}-\mathbf{y}\vert ^{2H+1 -\beta } \right) ^{-\frac{1}{2}}dtd\mathbf{x} ds d\mathbf{y} $$
and this is finite for $\beta >2H-1$. Thus almost surely the local time of  $u$ exists and is square integrable. \qed

\begin{remark}
It follows as a  consequence of Lemma 8.1 in \cite{Xiao1} that the local time of the solution $u$ admits the following $L^{2}$ representation
$$ L([a,b] \times [A,B], x)=\frac{1}{2\pi} \int_{\mathbb{R}} dz e^{-izx} \int_{[a,b] \times [A, B]} ds d\mathbf{y} e^{iu (s, \mathbf{y} ) z} $$
for every $x \in \mathbb{R}$.
\end{remark}

\subsection{Existence of the joint density for the solution in the Rosenblatt case}

It is possible to obtain the existence of the joint density of the random vector $ \left( u(t,\mathbf{x}), u(s,\mathbf{y}) \right) $ with $s\not=t$ or $\mathbf{x}\not= \mathbf{y}$  in the case when the wave equation (\ref{wave}) is driven by a Hermite sheet of order $q=2$ (the Rosenblatt sheet). The result is based on a criterium for the existence of densities for vectors of multiple integrals which has recently  been proven in \cite{NNP}.

Let us state our result.

\begin{prop}
Let $u(t,\mathbf{x}), t\geq 0, \mathbf{x} \in \mathbb{R} ^{d}$  be the mild solution to (\ref{wave}). Then for every $(t, \mathbf{x}) \not= (s, \mathbf{y})$, $(t, \mathbf{x}), (s, \mathbf{y} ) \in (0, \infty)\times \mathbb{R} ^{d}$, the random vector
$$\left( u(t,\mathbf{x}),  u(s, \mathbf{y}) \right) $$
admits a density.

\end{prop}
{\bf Proof: } Note that for every $t\geq 0$ and $\mathbf{x} \in \mathbb{R} ^{d}$, the random variable $u(t, \mathbf{x})$ is a multiple integral of order 2 with respect to the $d$-parametric Brownian sheet. A result present in \cite{NNP} states that a two-dimensional vector of multiple integrals of order 2 admits a density if and only if the determinant of the covariance matrix is strictly positive. Denote by $C(t,s, \mathbf{x}, \mathbf{y})$ the covariance matrix of $\left( u(t,\mathbf{x}),  u(s, \mathbf{y}) \right) $. The determinant of this matrix is the same for every $q\geq 1$, from the covariance structure of the Hermite sheet. It is clear that for $q=1$ obviously $\det C (t,s, \mathbf{x}, \mathbf{y})$ is strictly positive, since the vector $\left( u(t,\mathbf{x}),  u(s, \mathbf{y}) \right) $ is a Gaussian vector and hence admits a density when $(t, \mathbf{x}) \not= (s, \mathbf{y})$. This implies that   $\det C (t,s, \mathbf{x}, \mathbf{y})$  is also strictly positive for $q=2$ and so the vector $\left( u(t,\mathbf{x}),  u(s, \mathbf{y}) \right) $ admits a density also for $q=2$.
\qed


\begin{thebibliography}{99}


%



%


\bibitem{Abry1}
{Albin I.M.P. (1998):}
{\em A note on the Rosenblatt distributions.}
{ Statistics and Probability Letters, 40(1), 83-91.}

\bibitem{Abry2}
{Albin I.M.P. (1998): }
{\em On extremal theory for self similar processes.}
{Annals of Probability, 26(2), 743-793. }




\bibitem{ALP2002}
Ayache A., Leger S. \& Pontier M. (2002):
{\em Drap brownien fractionnaire. (French) [The fractional Brownian sheet].}
Potential Anal. {\bf 17}, no. 1, 3143.




\bibitem{BT}
Balan R.M. \& Tudor C. A. (2010):
{\em The stochastic wave equation with fractional noise: A random field approach.}
Stoch. Proc. Appl. {\bf 120}, 2468-2494.

\bibitem{BT10}
Bardet J.-M. \& Tudor C. A. (2010):
{\em A wavelet analysis of the Rosenblatt process: chaos expansion and estimation of the self-similarity parameter}. Stochastic Process. Appl. {\bf 120} 12, 2331-2362.

\bibitem{B1}
{Berman S. M. (1973):}
{\em Local nondeterminism and local times of Gaussian processes.}
Indiana. Univ. Math. J,  23, pp. 69-94.

\bibitem{B2}
{Berman S. M. (1969):}
{\em Local times and sample function properties of stationary Gaussian processes.}
Trans. Amer. Math. Soc. 137, 277 - 299.






\bibitem{Breton11}
Breton J.-C. (2011):
{\em On the rate of convergence in non-central asymptotics of the Hermite variations of fractional Brownian sheet.}
Prob. Math. Stats. {\bf 31}(2), 301-311.

%

\bibitem{Cencov}
{Cencov N.N. (1956):}
{\em Wiener random fields depending on several parameters.} Dokl. Akad. Nauk SSSR 106, 607-609.


\bibitem{CT12}
Clarke De la Cerda J. \& Tudor C. A. (2012):
{\em Hitting times for the stochastic wave equation with fractional-colored noise}.
Revista Matematica Iberoamericana, to appear, http://arxiv.org/abs/1203.3921v1.


\bibitem{CTV09}
Chronopoulou A., Tudor C. A. \& Viens F. G. (2009):
{\em Variations and Hurst index estimation for a Rosenblatt process using longer filters}.
Electron. J. Stat. {\bf 3}, 1393-1435.

\bibitem{CTV11-Hp}
Chronopoulou A., Tudor C.A. \& Viens F. G. (2011):
{\em Self-similarity parameter estimation and reproduction property for non-Gaussian Hermite processes}.
Commun. Stoch. Anal. {\bf 5} 1, 161-185.

%

%
%
\bibitem{DaSa1}
Dalang R. \& Sanz-Sol\'e M. (2010):
{\em Criteria for hitting probabilities with applications to systems of stochastic wave equations.}
Bernoulli {\bf 16}(4), 1343-1368.





\bibitem{DoMa}
{Dobrushin R. L. \& Major P. (1979):}
{\em Non-central limit theorems for non-linear functionals of Gaussian fields.}
{Zeitschrift f\"ur Wahrscheinlichkeitstheorie und verwandte Gebiete, {\bf 50}, 27-52. }







\bibitem{KP2009}
Kim Y. T. \& Park, H. S. (2009):
{\em Stratonovich calculus with respect to fractional Brownian sheet.}
Stoch. Anal. Appl. {\bf 27}, no. 5, 962-983.




\bibitem{MT07}
Maejima M. \& Tudor C. A. (2007):
{\em Wiener Integrals with respect to the Hermite process and a Non-Central Limit Theorem.}
 Stoch. Anal. Appl. {\bf 25} 5, 1043-1056.

\bibitem{NNP}
{Nourdin I., Nualart D. \& Poly G. (2012): }
{\em Absolute continuity and convergence of densities for random vectors on Wiener chaos.} arXiv:1207.5115.



\bibitem{N}
{Nualart D. (1995):}
{\em Malliavin Calculus and Related Topics.}
{Springer.}







\bibitem{PiTa1}
{Pipiras V. \& Taqqu M. S. (2001):}
{\em Integration questions related to the fractional Brownian motion.}
{Probability Theory and Related Fields, {\bf 118}, 251-281. }

\bibitem{PT10-Hp}
Pipiras V. \& Taqqu M. S. (2010):
{\em Regularization and integral representations of Hermite processes}.
Statist. Probab. Lett. {\bf 80} 23-24, 2014-2023.


\bibitem{RST10}
R\'eveillac A., Stauch M. \& Tudor C. A. (2011):
{\em Hermite variations of the fractional Brownian sheet.}
Stochastics and Dynamics,  12 (3), 21pp.

\bibitem{Rose}
{Rosenblatt M. (1960): }
{\em Independence and dependence. }
{Proceedings of the 4th Berkeley Symposium on Mathematical Statistics}, Vol. II, 1960, 431-443.


\bibitem{ST2008}
Sottinen T. \& Tudor C. A. (2008):
{\em Parameter estimation for stochastic equations with additive fractional Brownian sheet.}
Stat. Inference Stoch. Process. {\bf 11}, no. 3, 221-236.

\bibitem{Taqqu75}
{Taqqu M. S. (1975): }
{\em Weak convergence to the fractional
Brownian motion and to the Rosenblatt process. }
{Zeitschrift f\"ur Wahrscheinlichkeitstheorie und verwandte Gebiete. {\bf 31},  287-302. }


\bibitem{Taqqu79}
{Taqqu M. S. (1979):}
{\em Convergence of integrated processes of
arbitrary Hermite rank. }
{Zeitschrift f\"ur Wahrscheinlichkeitstheorie und verwandte Gebiete, {\bf 50} 53-83. }

\bibitem{treves75}
Treves F. (1975):
{\em Basic Linear Partial Differential Equations}.
Academic Press, New York.


\bibitem{Tudor08}
Tudor C. A. (2008):
{\em Analysis of the Rosenblatt process}.
ESAIM Probab. Stat. {\bf 12}, 230-257.

\bibitem{TV1}
{Tudor C. A. \& Viens F. G. (2003):}
{\em It\^o Formula and Local Time for the Fractional Brownian Sheet.}
Electronic Journal of Probability. {\bf 8}(14), 1-31.

\bibitem{Xiao1}
{Xiao Y. (2009): }
{\em Sample path  properties of anisotropic Gaussian random fields: A minicourse on stochastic partial differential equations.} Lecture Notes in Math. 1962, 145-212, Springer, Berlin.

\end{thebibliography}
 \end{document}